\documentclass[11pt]{article}
\usepackage{subeqn}
\usepackage         {amssymb}
\usepackage[dvips]  {graphicx}
\usepackage{color}

\def\mc{\mathcal}
\def\ee{\begin{equation}}
\def\eee{\end{equation}}
\def\bqn{\begin{eqnarray*}}
\def\eqn{\end{eqnarray*}}
\def\bnl{\begin{eqnarray}}
\def\enl{\end{eqnarray}}
\def\bma{\begin{bmatrix}}
\def\ema{\end{bmatrix}}
\def\bmx{\begin{matrix}}
\def\emx{\end{matrix}}
\def\ben{\begin{enumerate}}
\def\een{\end{enumerate}}
\def\bit{\begin{itemize}}
\def\eit{\end{itemize}}
\def\bei{\begin{itemize}}
\def\eei{\end{itemize}}
\def\bet{\begin{tabular}}
\def\eet{\end{tabular}}
\newcommand{\vn}[1]{\left|\left|#1\right|\right|}

\def \RR {\mathbb{R}}

\def\qedp{\hfill{$\blacksquare$}}
\def\qed{\hfill {$\square$}}
\def\salt{\vskip 0.4 true cm}
\newtheorem{theorem}{Theorem}

\newtheorem{definition}{Definition}

\newtheorem{remark}{Remark}

\newtheorem{lemma}{Lemma}

\title{\bf Synchronization in Networks of Identical Linear Systems}
\author{Luca Scardovi\thanks{Luca Scardovi is with the Department of Mechanical and Aerospace Engineering, Princeton University, USA. {scardovi@princeton.edu}. The work is supported in part by ONR grants N00014--02--1--0826 and N00014--04--1--0534.}\,\, and Rodolphe Sepulchre\thanks{Rodolphe Sepulchre is with the Department of Electrical Engineering and Computer Science,
University of Li\`ege, Belgium. {r.sepulchre@ulg.ac.be}. This paper presents research results of the Belgian Network DYSCO (Dynamical Systems, Control, and Optimization), funded by the Interuniversity Attraction Poles Programme, initiated by the Belgian State,  Science Policy Office. The scientific responsibility rests with its authors.
} 
}
\begin{document}
\maketitle
\begin{abstract}
The paper investigates the synchronization of a network of identical linear state-space models under a possibly time-varying and directed interconnection structure.
The main result is the construction of a  dynamic output feedback coupling that achieves synchronization if the decoupled systems have no exponentially unstable mode and if the communication graph is uniformly connected. The result can be interpreted as a generalization of classical consensus algorithms. Stronger conditions are shown to be sufficient -- but to some extent, also necessary -- to ensure synchronization with the diffusive static output coupling often considered in the literature.
\end{abstract}

\section{Introduction}
In these last years, consensus, coordination and synchronization problems have been popular subjects in systems and control, motivated by many applications in physics, biology, and engineering. These problems arise in multi-agent systems with the collective objective of reaching agreement about some variables of interest.

In the {\it consensus} literature, the emphasis is on the communication constraints rather than on the individual dynamics: the agents exchange information according to a communication graph that is not necessarily complete, nor even symmetric or time-invariant, but, in the absence of communication, the agreement variables usually have no dynamics. It is the exchange of information only that determines the time-evolution of the variables, aiming at asymptotic synchronization to a common value. The convergence of such consensus algorithms has attracted much attention in the recent years. It only requires a weak form of connectivity for the communication graph \cite{Mo,Mocdc,BlHeOlTs,JaLiMo,OlMu}.
 
In the {\it synchronization} literature, the emphasis is on the individual dynamics rather than on the communication limitations: the communication graph is often assumed to be complete (or all-to-all), but in the absence of communication, the time-evolution of the systems' variables can be oscillatory or even chaotic. The system dynamics can be modified through the information exchange, and, as in the consensus problem, the goal of the interconnection is to reach synchronization to a common solution of the individual dynamics \cite{Hale,PhSl,StSe,Po}.

Coordination problems encountered in the engineering world can often be rephrased as consensus or synchronization problems in which both the individual dynamics and the limited communication aspects play an important role. Designing interconnection control laws that can ensure synchronization of relevant variables is therefore a control problem that has attracted quite some attention in the recent years \cite{ScSaSe,SePaLe_limited,ScLeSe,NaLe,SaSeLe}.  

The present paper deals with a fairly general solution of the synchronization problem in the linear case.  Assuming $N$ identical  individual agents dynamics each  described by the linear state-space model $(A,B,C)$, the main result is the construction of a dynamic output feedback controller that ensures exponential synchronization to a solution of the linear system $\dot x = Ax$ under the following assumptions: (i) $A$ has no exponentially unstable mode, (ii) $(A,B)$ is stabilizable and $(A,C)$ is detectable, and (iii) the communication graph is uniformly connected. The result can be interpreted as a generalization of classical consensus algorithms, studied recently, corresponding to the particular case $A=0$ \cite{Mo,Mocdc}. The generalization includes the non-trivial examples of synchronizing harmonic oscillators or chains of integrators. 

The proposed dynamic controller structure proposed in this paper differs from the static diffusive coupling often considered in the synchronization literature, which requires more stringent assumptions on the communication graph.  For instance, the results in the recent paper \cite{Tun08} assume a time-invariant topology. 
The paper also provides sufficient conditions for synchronization by static diffusive coupling and illustrates on simple examples that  synchronization may fail under diffusive coupling when the stronger assumptions on the communication graph are not satisfied.

The paper is organized as follows. In Section \ref{sec:not} the
notation used throughout the paper is summarized, some preliminary results are reviewed and the
synchronization problem is introduced and defined. In Section \ref{sec:linearstate} the linear case is studied when state coupling among the systems is allowed while in Section \ref{sec:linearout} the output coupling case is considered. In Section \ref{sec:exten} we extend the main results to discrete-time linear systems and to periodic time-varying linear systems. Finally, in Section \ref{sec:example}, two-dimensional examples are reported to illustrate the role of the proposed dynamic controller in situations where  static diffusive coupling fails to achieve synchronization.    
\section{Preliminaries} \label{sec:not}
\subsection{Notation and Terminology}
Throughout the paper we will use the following notation. Given $N$ vectors $x_1, x_2, \ldots, x_N$ we indicate with $x$ the
stacking of the vectors, i.e. $x = [x_1^T, x_2^T, \ldots, x_N^T]^T$. We denote with $I_N$ the diagonal matrix of dimension
$N\times N$ and we define $1_{N}\triangleq [1,1,\ldots,1]^{T} \in \RR^{N}$. Given two matrices $A$ and $B$
we denote their Kronecker product with $A \otimes B$.  For notational
convenience, we use the convention $\tilde A_N = I_N \otimes A$ and $\hat A_N = A \otimes I_N$. We recall some properties of the Kronecker
product that will be used throughout the paper
\begin{subequations}
\label{kron}
\begin{eqnarray}
&&\hspace{-0.5 cm} A \otimes B \otimes C = A \otimes (B \otimes C) = (A \otimes B)\otimes C\label{kron1}\\
&&\hspace{-0.5 cm} A \otimes (B+C) = A \otimes B + A \otimes C\label{kron2}\\
&&\hspace{-0.5 cm}AB \otimes CD = (A \otimes C)(B \otimes D)\label{kron3}\\
&&\hspace{-0.5 cm}A\otimes B = (A \otimes I_p)(I_n \otimes B) = (I_m \otimes B)(A \otimes I_q)\label{kron4}\\
&&\hspace{-0.5 cm}AB \otimes I_n = (A\otimes I_n)(B \otimes I_n)\label{kron5}
\end{eqnarray}
\end{subequations}
where $\,A \in M^{mn}, \; B \in M^{pq}$.

\subsection{Communication Graphs}
Given a set of interconnected systems the communication topology is encoded through a  \emph{communication
graph}.  The convention is that system $j$ receives
information from system $i$ iff  there is a directed link from node $j$ to node $i$ in the communication graph.
Let ${\mc G}(t)=({\mc V},{\mc E}(t),A_d(t))$ be a time-varying weighted digraph (directed graph) where ${\mc V}=\{v_1,\ldots,v_N\}$ is
the set of nodes, ${\mc E}(t)\subseteq {\mc V} \times {\mc V}$ is the set of edges, and ${A}_{d}(t)$ is a weighted adjacency matrix
with nonnegative elements $a_{kj}(t)$. In the following we assume that $A_{d}(t)$ is piece-wise continuous and bounded and $a_{kj}(t)\in \{0\}\cup[\eta,\gamma],
\forall\, k, j,\,$ for some finite scalars $0<\eta\leq\gamma$ and for all $t\geq 0$. Furthermore $\{v_{k},v_{j}\} \in {\mc E}(t)$ if and only if $a_{kj}(t)\geq\eta$.   
The set of neighbors of node $v_k$ at time
$t$ is denoted by ${\mc N}_k(t) \triangleq \{v_j \in {\mc V} : a_{kj}(t)\geq\eta\}$. A path is a sequence of vertices such that for each of its vertices $v_{k}$ the next vertex in the sequence is a neighbor of $v_{k}$. Assume that there are no self-cycles i.e. $a_{kk}(t)=0, \,k=1,2,\ldots,N$, and for any
$t$.\\
The Laplacian matrix $L(t)$ associated to the graph ${\mc G}(t)$ is defined as
\[
l_{kj}(t)=\left\{
\begin{array}{ll}
\displaystyle \sum_{i=1}^{N} a_{ki}(t), &   j=k\\
-a_{kj}(t),       &   j\ne k.
\end{array}
\right.
\]
The in-degree (respectively out-degree) of node $v_k$ is defined as $d_k^{in}=\sum_{j=1}^N a_{kj}$ (respectively
$d_k^{out}=\sum_{j=1}^N a_{jk}$). The digraph ${\mc G}(t)$ is said to be {\em balanced} at time $t$ if the in-degree and the out-degree of each node
are equal, that is,
\[
\sum_{j=1}^{N} a_{kj} = \sum_{j=1}^{N} a_{jk}, \quad \quad k=1,\ldots,N.
\]
Balanced  graphs have the particular property that the symmetric part of their Laplacian matrix is nonnegative:  $L+L^T \ge 0$ \cite{Wi}.
We recall some definitions that characterize
the concept of connectivity for time-varying graphs.
\begin{definition}
The digraph ${\mc G}(t)$ is {\em connected} at time $t$ if there exists a node $v_{k}$ such that all the other nodes of the graph are connected to $v_{k}$ via a
path that follows the direction of the edges of the digraph.
\end{definition}
\begin{definition}
Consider a graph ${\mc G}(t)$. A node $v_k$ is said to be connected to node $v_j$ ($v_j \ne v_i$) in the
interval $I=[t_a,t_b]$ if there is a path from $v_k$ to $v_j$ which respects the orientation of the edges for the directed graph
$({\mc V},\cup_{t \in I} {\mc E}(t),\int_I A_{d}(\tau)d \tau)$.
\end{definition}
\begin{definition}
${\mc G}(t)$ is said to be uniformly connected if there exists a time horizon $T>0$ and an index $k$ such that for all $t$ all the
nodes $v_j$ ($j\neq k$) are connected to node $v_k$ across $[t, t+T]$.
\end{definition}
\subsection{Convergence of consensus algorithms}
Consider $N$ agents exchanging information about their state vector $x_{k}$, $k=1,\ldots,N$, according to a communication graph $\mc G(t)$. A classical consensus protocol in continuous-time is 
\ee \label{cons1} \dot{x}_k=\sum_{j=1}^{N} a_{kj}(t)( x_j- x_k), \quad \quad  k=1,\ldots,N.
\eee 

In discrete-time the analogous dynamics write
\ee \label{cons1d}
x_{k}(t+1) =x_{k}(t) -\epsilon_{k}\sum_{j=1}^{N}l_{kj}(t)x_{j}(t), \quad \quad k=1,\ldots,N 
 \eee
where $\epsilon_{k}\in (0,1/d_{k}^{in})$. 
Using the Laplacian matrix, (\ref{cons1}) and (\ref{cons1d}) can be equivalently expressed as \ee \label{cons2} \dot{
x}=-\hat{L}_n(t)\, x. \eee
and
 \ee \label{cons2d} 
x(t+1)=\left(I_{nN}-\hat \epsilon\, \hat{L}_n(t)\right)\, x(t), \eee
where $\hat \epsilon = \epsilon \otimes I_{n}$ and $\epsilon = \mbox{diag}(\epsilon_{1},\epsilon_{2}, \ldots,\epsilon_{N})$. 

Algorithms (\ref{cons2}) and (\ref{cons2d}) have
been widely studied in the literature and asymptotic convergence to a consensus value holds under mild assumptions on the communication topology. 
The following theorem summarizes the main results in \cite{Mo} and \cite{Mocdc}.
\begin{theorem}\label{th:ConsEucl}
Let $x_{k},\, k=1,2,\ldots,N$, belong to a finite-dimensional Euclidean space $W$. Let ${\mc G}(t)$ be a uniformly connected digraph and $L(t)$ the corresponding
Laplacian matrix bounded and piecewise continuous in time. Then the equilibrium sets of consensus states of (\ref{cons2}) and (\ref{cons2d}) are 
uniformly exponentially stable. 
Furthermore the solutions of (\ref{cons2}) and and (\ref{cons2d}) asymptotically converge to a consensus value $1_{N} \otimes \beta $ for some $\beta \in W$.
\qed
\end{theorem}
A general proof for Theorem \ref{th:ConsEucl} is based on the property that the convex hull of vectors $x_k \in W$ is non expanding along the solutions. 

\subsection{The Synchronization Problem} \label{sec:main}
Consider $N$ identical dynamical systems
\begin{subequations}\label{model1}
\begin{eqnarray}
\dot x_k &=& f(t,x_k,u_k)\\
y_k &=& h(x_k), 
\end{eqnarray}
\end{subequations}
for $k=1,2,\ldots,N$, where $x_k \in \RR^n$ is the state of the system, $u_k \in \RR^m$ is the control and $y_k\in \RR^p$ is the output.  We assume that the coupling among the systems involves only the output differences $y_k - y_j$ and the controller state differences $\, \xi_k-\xi_j$. According to the graph-theoretic description of the communication topology, two systems are coupled at time $t$ if there exists an edge connecting them in the associated (time-varying) communication graph ${\mc G}(t)$ at time $t$. We will call a control law \emph{dynamic} if it depends on an internal (controller) state, otherwise it is called \emph{static}. For the systems to be synchronized, the control action (that will depend on the coupling) must vanish asymptotically and must force the solutions of the closed-loop systems to asymptotically converge to a common solution of the individual systems. This leads to the formulation of the following problem:   

\emph{Synchronization Problem:}
Given $N$ identical systems described by the model (\ref{model1}) and a communication graph ${\mc G}(t)$, find a (distributed) control law such that the solutions of (\ref{model1}) asymptotically synchronize to a solution of the open-loop system $\,\dot x_0 =  f(t,x_0,0)$.
\qed

In the present paper we focus the attention on synchronization of linear time-invariant systems. Generalizations will be the subject of future work.   
\section{Synchronization of linear systems with state feedback}\label{sec:linearstate}
Consider   $N$ identical linear systems, each described by the linear model
\begin{equation}\label{s1}
\dot{x}_k = A x_k + B u_k, \quad \quad  k=1,2,\ldots,N,
\end{equation}
where $x_k \in \RR^n$ is the state of the system and $u_k \in \RR^m$ is the control vector. 
For notational convenience it is possible to rewrite (\ref{s1}) in compact form as
\[
\begin{array}{rcl}
\dot{x} &=& \tilde{A}_N x + \tilde{B}_N u.
\end{array}
\]
Theorem \ref{th:ConsEucl} can be interpreted as  a synchronization result for linear systems with $A=0$ and $B=I$. A straightforward generalization is as follows.   \salt 
\begin{lemma}\label{le:ful} 
Consider the linear systems (\ref{s1}). Let $B$ be a $n \times n$ nonsingular matrix and assume that all the eigenvalues of $A$ belong to the imaginary axis. Assume that the communication graph $G(t)$ is uniformly connected and the corresponding Laplacian matrix $L(t)$ piecewise continuous and bounded. Then the control law
\[
u_k = B^{{-1}}\sum_{j=1}^N a_{kj}(t) (x_j - x_k),\quad k=1,2,\ldots,N,
\]
exponentially synchronizes all the solutions of (\ref{s1}) to a solution of the system $\dot {x}_0 = A x_0  $. \qed
\end{lemma}\salt
\emph{Proof:} Consider the closed-loop system 
\[
\dot{x}_k = A x_k +\sum_{j=1}^N a_{kj}(t) (x_j - x_k), \quad \quad k=1,2,\ldots,N.\]
The change of variable
\[
z_k = e^{-A(t-t_{0})} x_k, \quad \quad k=1,2,\ldots,N,
\]
leads to 
\[
\dot{z}_k=-  A  e^{-A(t-t_{0})} x_k  +  e^{-A(t-t_{0})} A x_k +  e^{-A(t-t_{0})}\sum_{j=1}^N a_{kj}(t) (x_j - x_k) =  \sum_{j=1}^N a_{kj}(t) (z_j - z_k)
\]
or, in a compact form, 
\[
\dot{z} = - \hat L_n(t)\, z.
\]
From Theorem \ref{th:ConsEucl} the solutions $z_k(t),\, k=1,2,\ldots,N$, exponentially converge to a common value $x_0 \in \RR^n$
as $t\rightarrow \infty$, that is, there exist constants $\delta_1>0$ and $\delta_2>0$ such that for all $t_0$,
\[
\hspace{-0.032 cm}\vn{z_k(t)-x_0} \leq \delta_{1} e^{-\delta_{2} (t-t_{0})}\vn{z_k(t_0)-x_0},\quad \quad \forall t>t_{0}. 
\]
In the original coordinates, this means
\[
\vn{x_k(t)-e^{A(t-t_0)}x_0} \leq  \delta_{1} \vn{e^{A(t-t_0)}} e^{-\delta_{2} (t-t_{0})}\vn{x_k(t_0)-x_0},
\]
for every $t>t_{0}$.
Because all the eigenvalues of the matrix $A$ lie on the imaginary axis, there exists a constant  
$\delta_3 >0$ such that   
\[
\vn{x_k(t)-e^{A(t-t_0)}x_0} \leq \delta_{1}  e^{-\delta_{3} (t-t_{0})}\vn{x_k(t_0)-x_0}, 
\]
for every $t>t_{0}$, which proves that  all solutions exponentially synchronize to a solution of the open loop system. \qedp 
\begin{remark}
The result is of course unchanged if $A$ also possesses eigenvalues with a negative real part. Exponentially stable modes synchronize to zero, even in the absence of coupling. In contrast, the  situation of systems with some eigenvalues with a positive real part can be addressed is a similar way but it requires that the graph connectivity is sufficiently strong to dominate the instability of the system. This is clear from the last part of the proof of Lemma \ref{le:ful} where the exponential synchronization  in the $z$ coordinates must dominate the divergence of the unstable modes of $A$.                
\end{remark}
The assumption of a {\it square} (nonsingular) matrix $B$ in Lemma \ref{le:ful} is now weakened to a stabilizability assumption on the pair $(A,B)$. For an arbitrary stabilizing feedback matrix $K$,  consider the (dynamic) control law
\begin{equation}\label{cont1}
\begin{array}{rcl}
\dot \eta_{k} &=&\left(A + B K\right) \eta_{k} + \displaystyle \sum_{j=1}^{N}a_{kj}(t) \left(\eta_{j}-\eta_{k}+x_{k}-x_{j}\right), \\
u_{k}&=&  K \eta_{k},
\end{array}
\end{equation}
for $k=1,2,\ldots,N$, which leads to the closed-loop system
\begin{subequations}
\label{coup2}
\begin{eqnarray}
\dot{x} &=& \tilde A_N  x + \tilde B_N \tilde K_N \eta \label{coup2_1}\\
\dot \eta &=& \left(\tilde A_N +\tilde B_N \tilde K_N\right) \eta + \hat L_n(t) (x-\eta).\label{coup2_2}
\end{eqnarray}
\end{subequations}

\begin{theorem}\label{th:main}
Consider the system (\ref{s1}). Assume that all the eigenvalues of $A$ belong to the closed left-half complex plane. Assume that the pair $(A,B)$ stabilizable and let $K$ be a stabilizing matrix such that $A+BK$ is Hurwitz. 
Assume that the graph is uniformly connected and that the Laplacian matrix is piecewise continuous and bounded. Then the solutions of (\ref{coup2}) exponentially synchronize to a solution of the open loop system $\dot x_0 =  A x_0$. 
\qed
\end{theorem}  
\emph{Proof:} With the change of variable $s_k =x_{k} - \eta_k$ we can rewrite (\ref{coup2_2}) as
\[
\dot{s} = \tilde A_N s - \hat L_n(t) s,
\]
and the   closed-loop dynamics write
\begin{subequations}\label{coup}
\begin{eqnarray}
\dot{x} &=& \left(\tilde A_N +\tilde B_N  \tilde K_N  \right)x + \tilde B_N  \tilde K_Ns \label{coup_1}\\
\dot{s} &=& \tilde A_N  s - \hat L_n(t) s. \label{coup_2}
\end{eqnarray}
\end{subequations}
Observe that the two systems (\ref{coup_1}) and (\ref{coup_2}) are decoupled. Since the assumptions of Lemma \ref{le:ful} are satisfied for the sub-system (\ref{coup_2}), its solutions exponentially  synchronize to a  solution of   $\dot s_0 =  A  s_0$.  The subsystem (\ref{coup2_2})  is therefore an exponentially stable system driven by  an
input $ \hat L_n(t)s(t)$  that exponentially converges to zero.  As a consequence, its solution $\eta(t)$ exponentially converges to zero, which implies that  the solutions of
(\ref{coup2_1}) exponentially synchronize to a solution of $\dot x_0= A x_0$. \qedp

\section{Synchronization of linear systems with output feedback}\label{sec:linearout}
Consider a group of $N$ identical linear systems described by the linear model
\begin{subequations}\label{s1o}
\begin{eqnarray}
\dot{x}_k &=& A x_k + B u_k,\\
y_k &=& C x_k
\end{eqnarray}
\end{subequations}
for $k=1,2,\ldots,N$, where $x_k \in \RR^n$ is the state of the system,  $u_k \in \RR^m$ is the control vector, 
and $y_k \in \RR^p$ is the output. 

The state feedback controller of  Theorem \ref{th:main} is easily extended to an output feedback controller if we additionally assume that the pair $(A,C)$ is detectable. Pick an observer matrix $H$ such that $A+HC$ is Hurwitz and consider the output feedback controller
 \begin{subequations}\label{contobs}
\begin{eqnarray}
\dot \eta_{k} &=&\left(A +B K\right) \eta_{k} +\displaystyle \sum_{j=1}^{N}a_{kj}(t) \left(\eta_{j}-\eta_{k}+\hat x_{k}-\hat x_{j}\right) \\
\dot{\hat x}_{k}&=& A \hat x_{k} +B u_{k} +H (\hat y_{k}-y_{k})\\
u_{k}&=& K \eta_{k}\\
\hat y_{k} &=& C \hat x_{k},
\end{eqnarray}
\end{subequations}
for $k=1,2,\ldots,N$, where detectability is assumed and $H$ is a suitable observer matrix. 
The convergence analysis is similar to the one for Theorem \ref{th:main} and is mainly based on the observation that the estimation error is decoupled from the consensus dynamics. \salt
\begin{theorem}\label{th:mainout}
Assume that the system (\ref{s1o}) is stabilizable and detectable and that all the eigenvalues of $A$ belong to the closed left-half complex plane.
Assume that the communication graph is uniformly connected and the Laplacian matrix is piecewise continuous and bounded. Then for any gain matrices
$K$ and $H$ such that $A+BK$ and $A+HC$ are Hurwitz, the solutions of (\ref{s1o}) with the dynamic controller (\ref{contobs}) exponentially synchronize to a solution of $\,\dot x_0= A x_0$.\flushright \qed
\end{theorem}
\emph{Proof:}
Define $s_{k} = \hat x_{k} - \eta_{k}$ and $e_{k} = x_{k} - \hat x_{k}$, and rewrite the closed loop system as
\begin{eqnarray*}
\dot x &=& \left(\tilde A_{N} + \tilde B_{N} \tilde K_{N}\right)x  + \tilde B_{N}\tilde K_{N}\left(e + s\right)\\
\dot s&=& \tilde A_{N}s -\hat L_{n} s\\
\dot e&=& \left(\tilde A_{N} + \tilde H_{N} \tilde C_{N}\right) e.
\end{eqnarray*}
This system is the cascade of  the closed-loop system analyzed in the proof of Theorem \ref{th:main} with an exponentially stable estimation error dynamics, which proves the result.
\qedp 

Theorem \ref{th:mainout} provides a general synchronization result for linear systems but the solution requires a dynamic controller. For the sake of comparison, we provide a set of sufficient conditions to prove synchronization under a simple static output feedback (diffusive) interconnection. These sufficient conditions require stronger assumptions on the interconnection and assume a passivity property for the system $(A,B,C)$, that is, the existence of a symmetric positive definite matrix $P>0$ that verifies  
\begin{equation} \label{pass}
PA + A^TP \leq 0,\quad B^TP = C.
\end{equation}

Passity conditions have been considered previously  
in \cite{Arcak} (where it is assumed that the communication topology is bidirectional and strongly connected) and in \cite{StSe} (where
synchronization is studied for a class of (nonlinear) oscillators  assuming that  the
communication topology is time-invariant and balanced).  Assumptions A1 and A2 below lead to a time-varying extension of the results in \cite{StSe} and \cite{Arcak} in the special case of linear systems. \salt
 \begin{theorem}\label{th:output}
Consider system (\ref{s1o}) with the static output feedback control laws
\[
u_k = \sum_{j=1}^N a_{kj}(t) (y_j - y_k),\quad \quad k=1,2,\ldots,N.
\]
Let the graph Laplacian matrix $L(t)$ be piecewise continuous and bounded. Then exponential synchronization to a solution of $\dot x_0 = A x_0$  is achieved under either one of the following assumptions:   \\
A1. The system $(A,B,C)$  is passive and observable, the communication graph is connected and balanced at each time;\\
A2. The system $(A,B,C)$ is passive and observable, the communication graph is symmetric, i.e. the Laplacian matrix can be factorized as $L=DD^T(t)$, and the pair $(\tilde A_N, \hat D_p^T(t) \tilde C_N)$ is uniformly observable. \qed
\end{theorem}

\emph{Proof:} Supppose that assumption A1 holds and consider the matrix $P$ solution of (\ref{pass}).\\
 Consider the Lyapunov function
\begin{equation}\label{V1}
V(x) =  \frac{1}{2}(\hat \Pi_n x)^T \tilde P_N (\hat \Pi_n x),
\end{equation}
the derivative along the solutions of the closed loop system is
\[
\dot V(x) =\frac{1}{2} \dot x^T \hat \Pi_n  \tilde P_N \hat \Pi_n \tilde A_N  x + \frac{1}{2}x^T \hat \Pi_n  \tilde P_N \hat \Pi_n \tilde A_N  \dot x.
\]
By using the commutation property (\ref{kron4}) of Kronecker product and the passivity relation (\ref{pass}) we obtain
\begin{equation}\label{dotVprov}
\begin{array}{rcl}
\dot V(x) &=&\displaystyle \frac{1}{2}x^T \hat \Pi_n (\tilde P_N \tilde A_N +\tilde A_N^T \tilde P_N)\hat
\Pi_n  x - x^T \tilde C^T_N \Pi_p 
\hat L^{{\footnotesize \mbox{sym}}}_p(t) \hat \Pi_p y\\
&\leq&  -y^T\hat \Pi_p \hat L^{{\footnotesize \mbox{sym}}}_p(t) \hat \Pi_p y.
\end{array}
\end{equation}
Because the graph is balanced, the matrix $L^{\mbox{\footnotesize sym}}(t)\triangleq (L(t)+L^T(t)) /2$ is
positive semi-definite for each $t$ and 
\[
(\hat \Pi_p y)^T \hat L^{{\mbox{\footnotesize sym}}}_p(t) \hat \Pi_p y \geq \lambda^{*}_2 \vn{\hat \Pi_p y}^2,
\]
where $\lambda^{*}_2 = \inf_{t} \lambda_{2}(t)$, and $\lambda_{2}(t)$  is the  algebraic connectivity of the graph at time $t$. Note that $\lambda^{*}_2>0$ because the graph is connected at each time $t$ and the values of the adjacency matrix related to the connected components are assumed to be bounded away from zero (see Section \ref{sec:not}).  This allows to rewrite (\ref{dotVprov}) as
\begin{equation}\label{dotVfin}
\dot V(x) \leq - \lambda^{*}_2 \vn{\hat \Pi_p y}^2, \quad \lambda^{*}_2>0.  
\end{equation}
Integrating (\ref{dotVfin}) over the interval $[t_0, t_0 + T]$ where $T>0$ is arbitrary, we obtain
\[
\displaystyle \int_{t_0}^{t_0 + T} \dot V dt \leq -\lambda^{*}_2 \displaystyle \int_{t_0}^{t_0 + T} \vn{\hat \Pi_p y}^2  dt
\leq -\gamma \lambda^{*}_2 \vn{\hat \Pi_n x(t_0)}^2, \quad \gamma>0,
\]
for all $x(t_{0})$, where the last inequality follows from the observability condition of the pair $(A,C)$.
%
We conclude from a standard Lyapunov argument that the solutions exponentially
converge to the invariant subspace
 \begin{equation}\label{invariant}
 \left\{x\in R^{nN}:\; x_{k} = \frac{1}{N} \sum_{j=1}^{N}x_{j},\; k=1,2,\ldots,N\right\},
 \end{equation}
and therefore they exponentially synchronize. To prove that they actually synchronize to a solution of the open-loop system it is sufficient to observe that the coupling vanishes in (\ref{invariant}). This implies that the solutions converge to the  $\omega$-limit sets of the uncoupled system that belong to (\ref{invariant}), concluding the first part of the proof.       
 
Assume that assumption A2 holds. First observe that from the symmetry of the communication graph the Laplacian matrix can be
factorized as $L(t) = DD^T(t)$. Uniform observability of the pair $(\tilde A_N, \hat D_p^T
\tilde C_N)$ means that for all $t_{0}>0$ there exist positive constants $T$ and $\alpha$ (independent from $t_{0}$) such that
\[
\int_{t_{0}}^{t_{0} + T} \tilde \Phi_N(t,t_{0})^T \tilde C_N^T \hat D_p \hat D_p^T(\tau) \tilde C_N \tilde \Phi_N(t,t_{0})dt \geq
\alpha I_{nN},
\]
where $\Phi(t,\tau)$ is the transition matrix. This implies that the system
\begin{subequations}\label{sys_obs}
\begin{eqnarray}
\dot x &=& \tilde A_{N} x\\
z &=& \hat D_p^T(t) \tilde C_N x,
\end{eqnarray}
\end{subequations}
is uniformly observable. Applying output injection to system (\ref{sys_obs}) we obtain
\begin{eqnarray*}
\dot x &=& \tilde A_{N} x - K(t)\tilde D_p^T \tilde C_N x\\
z &=& \hat D_p^T(t) \tilde C_N x.
\end{eqnarray*}
Choose $K(t)\triangleq \tilde P_N^{-1}\tilde C_N^T \hat D_p^T(t)$ and observe that, since $L(t)$ is bounded, $K(t)$ belongs to
$L_2(t,t+T)$. Then output injection preserves observability (see \cite{AESePe} and references therein) and the system
\begin{eqnarray*}
\dot x &=& \tilde A x - \tilde B_N\hat D_p^T \hat D_p(t) \tilde C_N x\\
z &=& \hat D_p^T(t) \tilde C_N x
\end{eqnarray*}
is still uniformly observable (here we have also used the passivity condition $\tilde C_N=\tilde B_N^T \tilde P_N$). Therefore for all
$t_{0}>0$ there exist positive constants $T$ and $\beta$ (independent from $t_{0}$) such that for every $x(0)\neq 0$
\[
\int_{t_{0}}^{t_{0} + T}\vn{z}^2dt= \int_{t_{0}}^{t_{0} + T} y(t)^T \hat D_p \hat D_p^T(t) y(t) dt \geq \beta.
\]

Consider the Lyapunov function (\ref{V1}). Integrating its time derivative over the interval
$[t_0, t_0 + T]$ where $T>0$ is arbitrary we obtain
\[
\displaystyle  \int_{t_0}^{t_0 + T} \dot V dt \leq - \displaystyle  \int_{t_0}^{t_0 + T} \vn{\hat \Pi_p \hat D_{p} \hat D_{p}^T(t)  y}^2 dt
\leq -\sigma \vn{\hat \Pi_n x(t_{0})}^2, \quad \sigma>0.
\]
We conclude from standard Lyapunov results that the solutions
asymptotically synchronize. The rest of the proof is equivalent to the end of the proof under Assumption A1.
\qedp

\section{Extensions and Generalizations} \label{sec:exten}
In the previous sections the results have been presented for \emph{time-invariant} linear systems in \emph{continuous time}. For the sake of completeness, we briefly discuss straightforward extensions to discrete-time systems and periodic systems.  

\subsection{Discrete-Time Linear Systems}
The first step is to provide a discrete-time counterpart of Lemma \ref{le:ful}.
Consider the discrete-time linear system
\begin{equation}\label{disc1}
x_{k}(t+1) = A x_{k}(t) +Bu_{k},  \quad\quad t=1,2\ldots, \quad k=1,\ldots,N.
\end{equation}
From Theorem \ref{th:ConsEucl} we know that the solutions of the system
\[
z_{k}(t+1) =z_{k}(t) -\epsilon_{k}\sum_{j=1}^{N}l_{kj}(t)z_{j}(t), \quad \quad k=1,\ldots,N 
 \]
where $\epsilon_{k}\in (0,1/d_{k}^{in})$, asymptotically converge to a consensus value if the assumptions of Theorem  \ref{th:ConsEucl} are satisfied.
With the change of variable $x_{k} = A^{(t-t_{0})}z_{k},\, t>t_{0},$ we obtain
\begin{equation}\label{dimdisc}
x_{k}(t+1) = A^{(t+1-t_{0})}z_{k}(t+1) =   Ax_{k}(t) +\epsilon_{k}A\sum_{j=1}^{N}l_{kj}(t)x_{j}(t), \quad \quad k=1,\ldots,N. 
\end{equation}
Identifying (\ref{disc1}) and (\ref{dimdisc}) results in the control law 
\ee \label{contdisc}
u_{k} = \epsilon_{k}\, B^{-1}A\sum_{j=1}^{N}l_{kj}(t)x_{j}(t), 
\eee
where we assumed that $B$ is invertible.
Assume now that the eigenvalues of $A$ belong to the unit circle (in the complex plane). Then there exist $\gamma>0$ and $0<q<1$ such that for all $t_{0}$ 
 \[
 \vn{x_{k} - A^{(t-t_{0})}x_{0}} \leq \vn{A^{(t-t_{0})} }   \vn{z_{k} - x_{0}}\leq \gamma \, q^{(t-t_{0})}    \vn{x_{k}(t_{0}) - x_{0}},\quad \quad  k=1,\ldots, N,\quad t>t_{0}.  
  \]
This shows that the solutions of system (\ref{disc1}) equipped with (\ref{contdisc}) synchronize to a solution of the open-loop system
$x_{0}(t+1) = A x_{0}(t)$.
This result is summarized in the following theorem.
\begin{lemma}\label{le:disc1}
Consider the system (\ref{disc1}). Let $B$ be a $n \times n$ nonsingular matrix and assume that all the eigenvalues of $A$ belong to the boundary of the unitary closed disk (in the complex plane). Assume that the communication graph ${\mc G}(t)$ is uniformly connected and the corresponding Laplacian matrix $L(t)$ piecewise continuous and bounded. Then the control law
\[
u_{k} = \epsilon_{k}B^{-1}A\sum_{j=1}^{N}l_{kj}(t)x_{j}(t), \quad \quad k=1,2,\ldots,N,\quad \epsilon_{k}\in (0,1/d_{k}^{in}),
\] 
exponentially synchronizes all the solutions of (\ref{disc1}) to a solution of the system $x_{0}(t+1) = A x_{0}(t)$.  \qed 
\end{lemma}
Thanks to Lemma \ref{le:disc1}, we can recast Theorem \ref{th:main} and Theorem \ref{th:mainout} in a discrete-time setting. For the sake of compactness we only report the output-feedback case (the state feedback case is just a particular case that can be easily derived by the reader). 

Consider the system
\begin{subequations}\label{disc3}
\begin{eqnarray}
x_{k}(t+1) &=& A x_{k}(t) +Bu_{k}(t),\\
y_{k}(t)&=&C x_{k}(t),
\end{eqnarray}
\end{subequations}
for $k=1,2,\ldots,N$, and the discrete-time version of (\ref{contobs})
\begin{subequations}
\label{contobsd}
\begin{eqnarray}
\eta_{k}(t+1) &=& \left(A  + BK\right)\eta_{k}(t) +\epsilon_{k}A\sum_{j=1}^{N}l_{kj}(t) \left(\hat x_{j}(t)-\eta_{j}(t)\right),\\
\hat x_{k}(t+1) &=& A\hat x_{k}(t) +H\left(y_{k}(t)-\hat y_{k}(t)\right),
\end{eqnarray}
\end{subequations} 
for $k=1,2,\ldots,N$, where $\hat y_{k}(t) = C \hat x_{k}(t)$. 
\begin{theorem}\label{th:disc3}
Assume that the system (\ref{disc3}) is stabilizable and detectable and that all the eigenvalues of $A$ belong to the the closed unitary disk in the complex plane.
Assume that the communication graph ${\mc G}(t)$ is uniformly connected and the Laplacian matrix is piecewise continuous and bounded. Then for any gain matrices
$K$ and $H$ such that $A+BK$ and $A+HC$ are Schur matrices, the solutions of (\ref{disc3}) with the dynamic controller (\ref{contobsd}) exponentially synchronize to a solution of $x_0(t+1)= A x_0(t)$.\qed
\end{theorem} 
The proof of Theorem \ref{th:disc3} is straightforward adaptation of the continuous-time counterpart and is therefore omitted. 

\subsection{Periodic Linear systems}
Periodic linear systems, naturally arise in a number of contexts in engineering, physics, and biology \cite{Periodic}.  Periodic models are
of large interest also in time-series analysis, economy and
finance, and in all other cases when seasonal phenomena has to
be taken in account.
Furthermore they arise when linearization of a nonlinear system along a periodic solution is analyzed. Therefore it is not difficult to figure out possible applications when synchronization of such models can be of interest.    
The results presented in this section follow from the well-known Floquet theory, where the properties of linear periodically time-varying systems are studied via a state transformation into a new coordinate system in which the system matrix becomes time-invariant. The eigenvalues of this matrix are called the \emph{characteristic exponents} of the original time-varying system matrix. In the following the continuous-time case is analyzed, the discrete-time case follows the same lines and is omitted for the sake of brevity.  
Consider the time-varying extension of (\ref{s1})
\begin{equation}\label{s1tv}
\dot{x}_k = A(t) x_k + B u_k, \quad \quad  k=1,2,\ldots,N,
\end{equation}
where $x_k \in \RR^n$ is the state of the system and $u_k \in \RR^m$ is the control vector. 
The following Theorem generalizes Lemma \ref{le:ful} to periodic linear systems. \salt
\begin{lemma}\label{le:per}
Consider the linear systems (\ref{s1tv}) where $A(\cdot)$ is periodic of period $T$ and let $B$ be a $n \times n$ nonsingular matrix. Assume that the \emph{characteristic exponents} of $A(\cdot)$ belong to the closed left half complex plane. Assume that the communication graph $\mc{G}(t)$ is uniformly connected and the corresponding Laplacian matrix $L(t)$ piecewise continuous and bounded. Then the control law
\[
u_k = B^{{-1}}\sum_{j=1}^N a_{kj}(t) (x_j - x_k),\quad k=1,2,\ldots,N,
\]
exponentially synchronizes all the solutions of (\ref{s1tv}) to a solution of the system $\dot {x}_0 = A(t) x_0  $. \qed
\end{lemma} \salt
\emph{Proof:} 
By using Floquet Theory (see for instance \cite{Perko}), there exists a time varying linear transformation
\[
s_{0}(t)=Q^{-1}(t)x_{0}(t)
\] 
where $Q(t)$ is continuous, non-singular and periodic of period $T$, such that the linear-time varying system 
\[
\dot x_{0} = A(t) x_{0}
\]
reduces to the linear time-invariant system
\[
\dot s_{0} = \Omega s_{0},
\]
where $\Omega$ is a constant matrix and its eigenvalues are the characteristic exponents of the original system. Moreover the transition matrix can be written as
\[
\Phi(t,t_{0}) = Q(t-t_{0})e^{\Omega (t-t_{0})}.
\]
At this point, following the same lines of the proof of Lemma \ref{le:ful}, we observe that with the linear transformation $z_{k} = \Phi(t_{0},t) x_{k}$ the solutions (in the $z$ coordinates)  asymptotically converge to a consensus state $x_{0}$.       
It follows that there exist constants $\delta_1>0$ and $\delta_2>0$ such that for all $t_0$
\[
\vn{x_k(t)- Q(t-t_{0})e^{\Omega(t-t_{0})}x_0} \leq  \delta_{1} \vn{e^{\Omega(t-t_{0})}}\vn{Q(t-t_{0})} e^{-\delta_{2} (t-t_{0})}\vn{x_k(t_0)-x_0},\quad \quad \forall t>t_{0}. 
\]
Since $Q(t)$ is continuous and periodic it is also bounded. Moreover the eigenvalues of $\Omega$ (the characteristic exponents of $ A(\cdot)\ $) are in the close left half complex plane and therefore there are no exponentially unstable modes.  We conclude that there exist constants  
$\delta_3 >0$ and $\delta_{4}>0$ such that for all $t_{0}$
\[
\vn{x_k(t)- Q(t-t_{0})e^{B(t-t_{0})}x_0} \leq  \delta_{3} e^{-\delta_{4} (t-t_{0})}\vn{x_k(t_0)-x_0},\quad \quad \forall t>t_{0},
\]
and therefore the solutions exponentially synchronize to a solution of the system $\dot x_{0} = A(t) x_{0}$. \qedp  

Following the same steps as in Section \ref{sec:linearstate}, it is straightforward to translate the results of Theorem \ref{th:main} and Theorem \ref{th:mainout} to the periodic case. We leave the extension to the interested reader. 
 
\section{Examples}\label{sec:example}
The conditions of Theorem \ref{th:output} are only {\it sufficient} conditions for exponential synchronization under diffusive coupling. We provide two simple examples to illustrate that these conditions are not far from being necessary when considering time-varying and directed graphs and that the {\it internal model} of the dynamic controller (\ref{cont1}) plays an important role in such situations.

\emph{Example 1: Synchronization of harmonic oscillators}\\
\begin{figure}[tb]
\centering
 \includegraphics[scale=0.34]{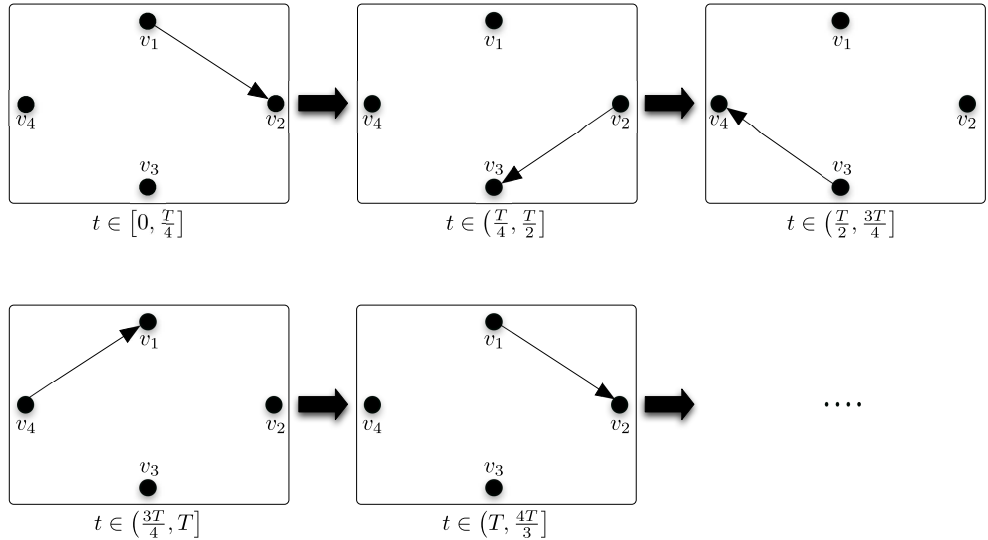}
\caption{The time-varying communication topology used in Example 1 and Example 2.} \label{fig:topol}
\end{figure}
Consider a group of $N$ harmonic oscillators
\begin{eqnarray*}
\dot x_{1k} &=& x_{2k}        \\
\dot x_{2k} &=& -x_{1k} + u_k,
\end{eqnarray*}
for $k=1,2,\ldots,N$, which corresponds to system (\ref{s1}) with 
$$ A = \left ( \begin{array}{cc} 0 & 1  \\ -1 & 0 \end{array} \right ), \;  B = \left ( \begin{array}{c} 0   \\ 1 \end{array}\right ). $$
 The assumptions of Theorem \ref{th:main} are satisfied: $A$  is Lyapunov stable and $(A,B)$ is stabilizable. Choosing the stabilizing gain $K= (0 \;- 1)$, the dynamic control law (\ref{cont1}) yields  the closed-loop system 
\begin{eqnarray*}
 \dot x_{1k} &=&x_{2k}       \\
\dot x_{2k}&=& -x_{1k} - \eta_{2k}  \\
\dot \eta_{1k} &=& \eta_{2k}  + \sum_{j=1}^N a_{kj}(t)(\eta_{1j}  - \eta_{1k} +  x_{1k} - x_{1j})  \\
\dot \eta_{2k} &=& -\eta_{1k} - \eta_{2k} + \sum_{j=1}^N a_{kj}(t)(\eta_{2j}  - \eta_{2k} +  x_{2k} - x_{2j}),
\end{eqnarray*}
for $k=1,2,\ldots,N$.
\begin{figure}[t]
\centering
\begin{tabular}{c c}
\hspace{-0.7 cm} \includegraphics[scale=0.35]{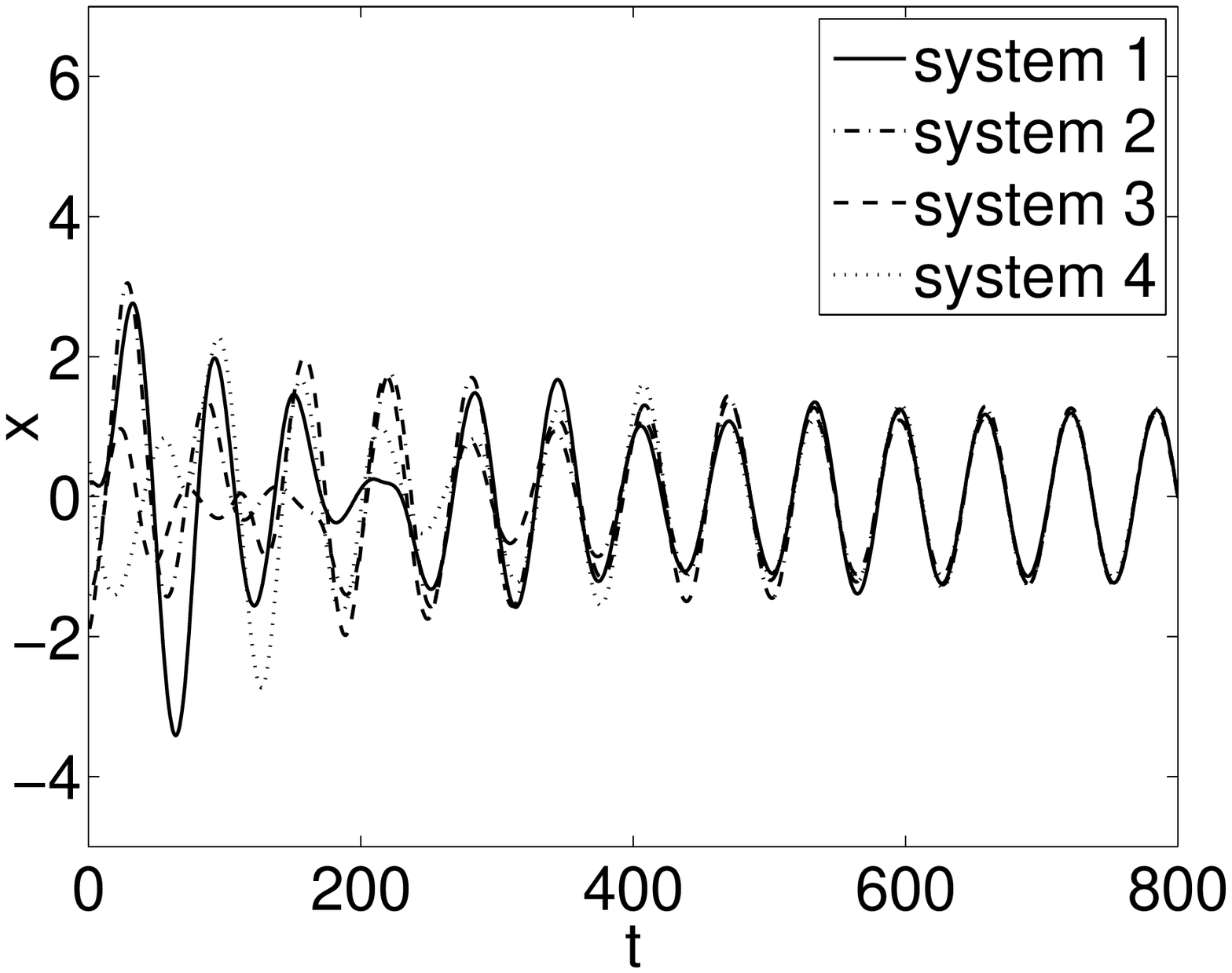}&\hspace{-0.4 cm}\includegraphics[scale=0.35]{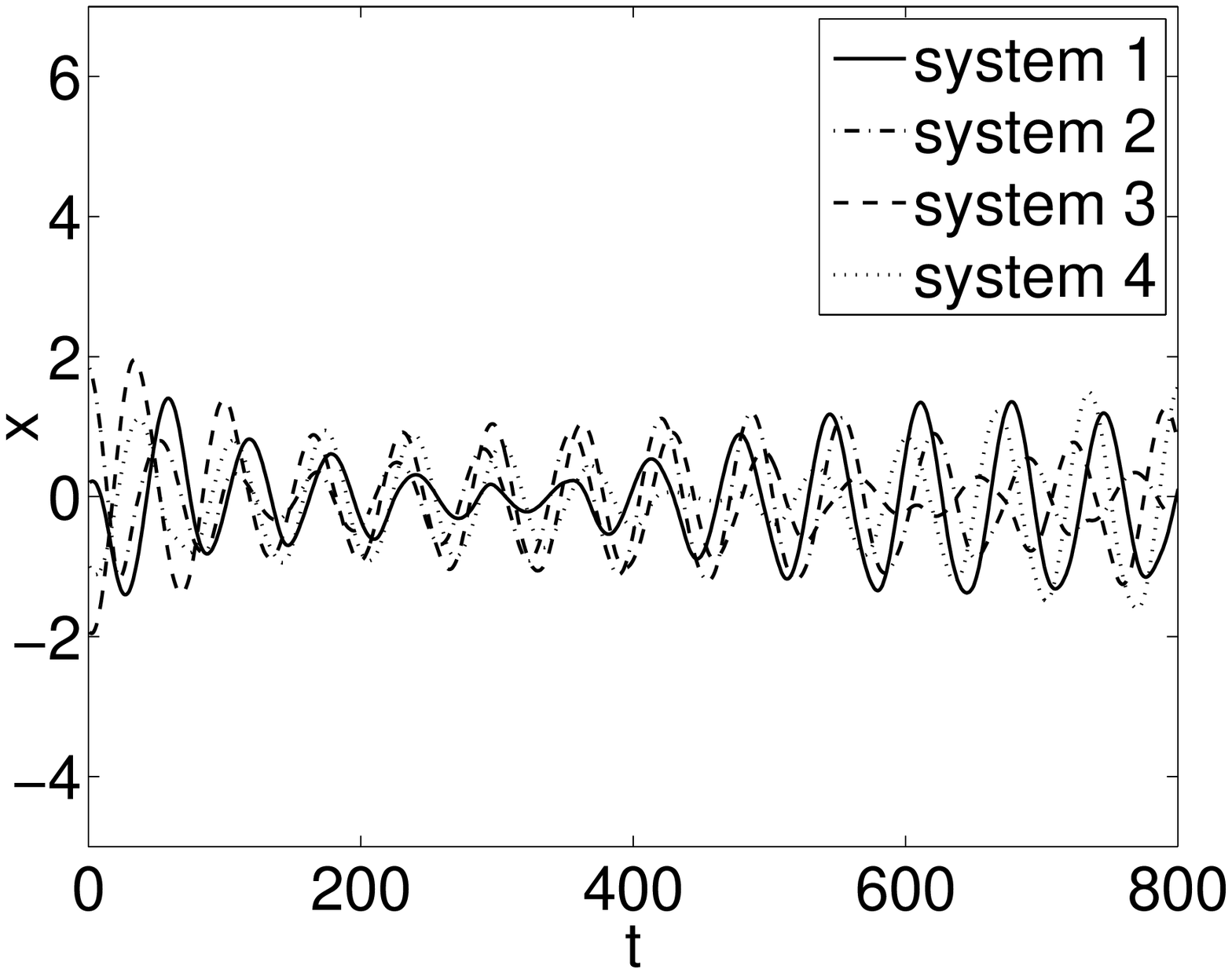}\\
\end{tabular}
\caption{First component of the solutions of the closed loop harmonic oscillators by using the dynamic control law (to the left) and the static control law (\ref{static}) (to the right). The dynamic control ensures exponential synchronization. In contrast, synchronization is not observed with the  diffusive interconnection.
} \label{fig:harm}
\end{figure}
Theorem \ref{th:main}  ensures exponential synchronization of the  oscillators  to a solution of the harmonic oscillator if the graph is uniformly connected.
Fig. \ref{fig:harm} illustrates the simulation of a group of $4$ oscillators coupled according to the time-varying communication topology shown in Fig.
\ref{fig:topol} (the period $T$ is set to $7$ sec). The dynamic control ensures exponential synchronization. In contrast, synchronization is not observed with the  diffusive interconnection
\begin{equation}\label{static}
u_k = \sum_{j=1}^N a_{kj}(t)(x_{2j} - x_{2k}).
\end{equation}
The system $(A,B,-K)$ is nevertheless passive, meaning that    stronger assumptions on the communication graph would ensure synchronization with the diffusive coupling (\ref{static}). We mention the recent result \cite{Tun08} that proves (in discrete-time) synchronization of harmonic oscillators with diffusive coupling under the assumption that the graph is time-invariant and connected.
The following example illustrates an analog scenario with unstable dynamics.

\emph{Example 2: Consensus for double integrators}\\
Consider a group of $N$ double integrators
\begin{eqnarray*}
\dot x_{1k} &=& x_{2k}        \\
\dot x_{2k} &=&   u_k,
\end{eqnarray*}
for $k=1,2,\ldots,N$, which corresponds to system (\ref{s1}) with 
$$ A = \left ( \begin{array}{cc} 0 & 1  \\ 0 & 0 \end{array} \right ), \;  B = \left ( \begin{array}{c} 0   \\ 1 \end{array} \right ). $$
 The assumptions of Theorem \ref{th:main} are satisfied: the two eigenvalues of $A$  are zero and $(A,B)$ is stabilizable. Choosing the stabilizing gain $K= (-1 \; -1)$, the dynamic control law (\ref{cont1}) yields  closed-loop system 
\begin{eqnarray*}
\dot x_{1k} &=&x_{2k}       \\
\dot x_{2k} &=&-\eta_{1k} - \eta_{2k}  \\
\dot \eta_{1k}&=&\eta_{2k}  + \sum_{j=1}^N a_{kj}(t)(\eta_{1j}  - \eta_{1k} +  x_{1k} - x_{1j})  \\
\dot \eta_{2k}&=& -\eta_{1k} - \eta_{2k} + \sum_{j=1}^N a_{kj}(t)(\eta_{2j}  - \eta_{2k} +  x_{2k} - x_{2j}),
\end{eqnarray*}
for $k=1,2,\ldots,N$.
\begin{figure}[t]
\centering
\begin{tabular}{c c}
 \hspace{-0.7 cm}  \includegraphics[scale=0.35]{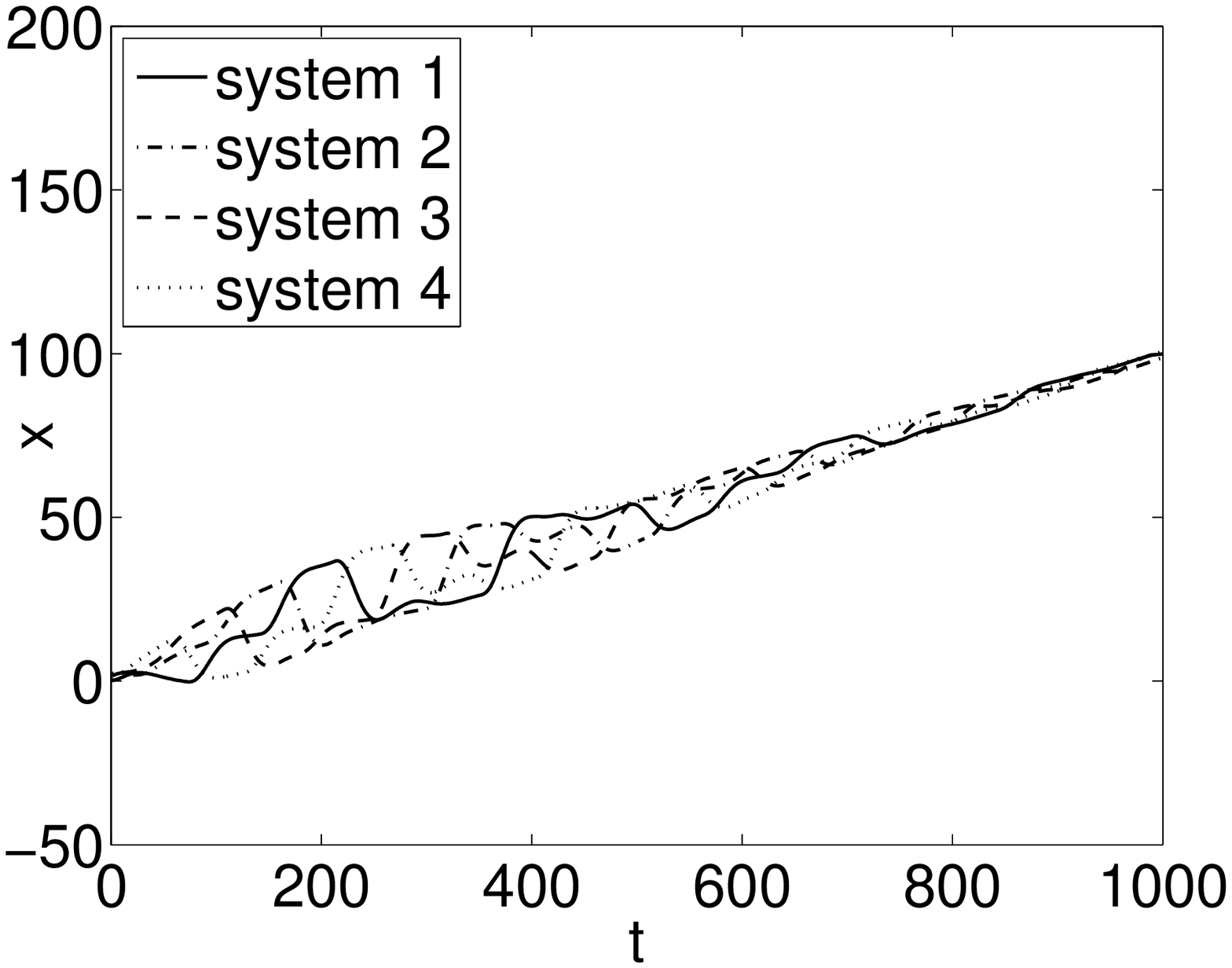}&\hspace{-0.4 cm}\includegraphics[scale=0.35]{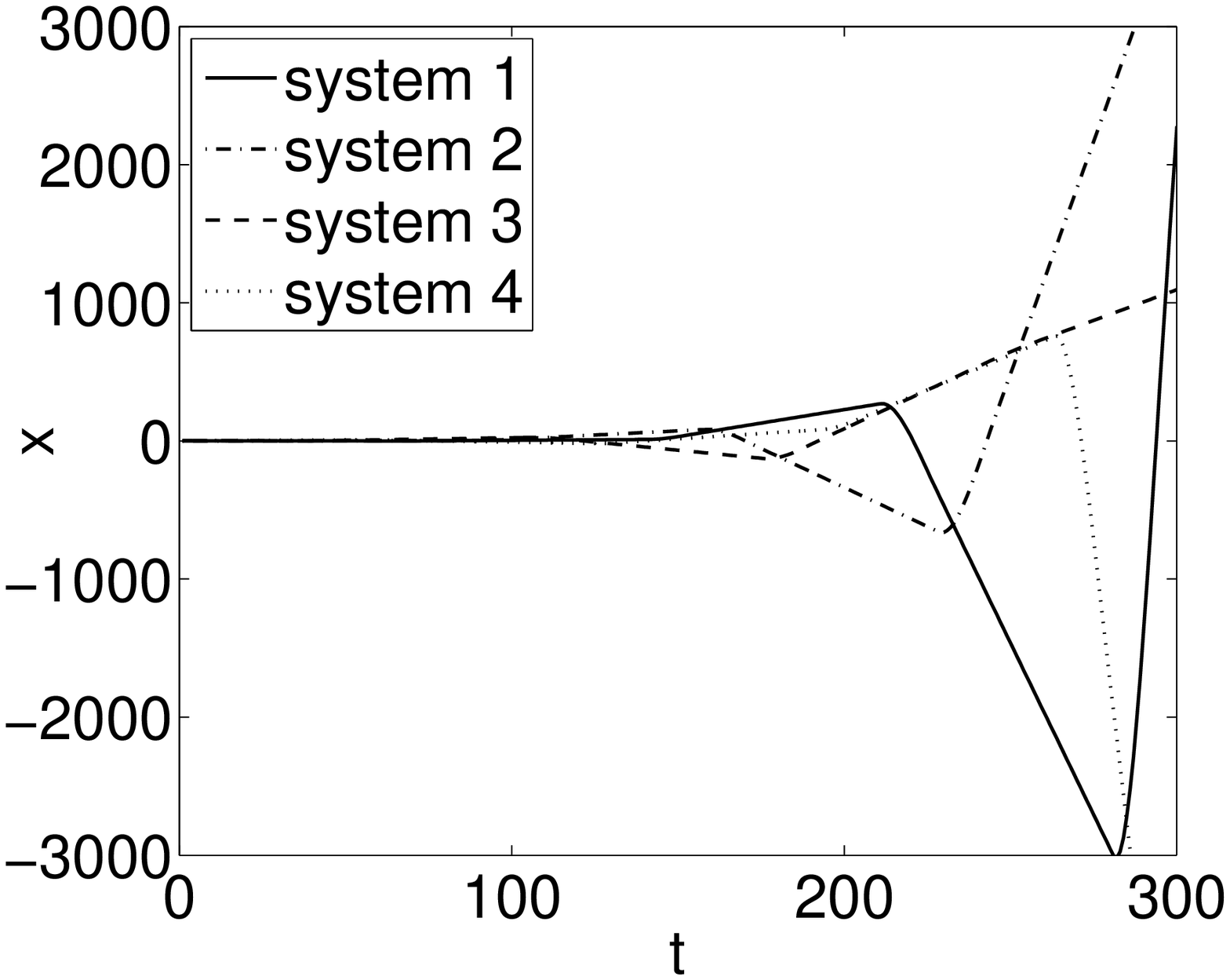}
\end{tabular}
\caption{First component of the solutions of the closed loop double integrators by using the dynamic control law (to the left) and the static control law
(\ref{staticdi}) (to the right).  The dynamic control ensures exponential synchronization. In contrast synchronization is not observed with the  diffusive interconnection.}
  \label{fig:double}
\end{figure}
Theorem \ref{th:main}  ensures exponential synchronization to a solution of the double integrator if the graph is uniformly connected.
Fig. \ref{fig:double} illustrates the simulation of a group of $4$ double integrators coupled according to the time-varying communication topology shown in Fig.
\ref{fig:topol} (the period $T$ is set to $2$ sec). The dynamic control ensures exponential synchronization. In contrast, synchronization is not observed with the  diffusive interconnection
\begin{equation}\label{staticdi}
u_k = \sum_{j=1}^N a_{kj}(t)(y_{j} - y_{k}), \quad \quad y_k = x_{1k}+x_{2k}.
\end{equation}
The matrix $A-\alpha BC$ is nevertheless stable for every $\alpha>0$, suggesting that stronger assumptions on the communication graph would ensure synchronization.

\section{Conclusion and future work}
In this paper the problem of synchronizing a network of identical linear systems described by the state-space model $(A,B,C)$ under general interconnection topologies has been addressed. A dynamic controller ensuring exponential convergence of the solutions to a synchronized solution of the decoupled systems is provided assuming that  (i) $A$ has no exponentially unstable mode, (ii) $(A,B)$ is stabilizable and $(A,C)$ is detectable, and (iii) the communication graph is uniformly connected. Stronger conditions are shown to be sufficient (and, to some extent, also necessary) to ensure synchronization with the often considered static diffusive output coupling.  The extension of the proposed technique for synchronization of nonlinear systems is the subject of ongoing work.

\end{document}